\input amstex 
\documentstyle{amsppt}
\input bull-ppt
\keyedby{bull326/jad}

\topmatter
\cvol{27}
\cvolyear{1992}
\cmonth{October}
\cyear{1992}
\cvolno{2}
\cpgs{292-297}
\define\Ext{\operatorname{Ext}}
\title Voiculescu Theorem, Sobolev Lemma, and\\Extensions 
of smooth algebras
\endtitle
\author Xiaolu Wang\endauthor
\dd the memory of Xian-Rong Wang\enddd
\shorttitle{Voiculescu theorem and Sobolev lemma}
\address Department of Mathematics, University of 
Maryland, College 
Park, Maryland  20742\endaddress
\ml xnw\@karen.umd.edu\endml
\date March 24, 1992\enddate
\comm{Xian-Rong Wang}
\abstract We present the analytic foundation of a unified 
B-D-F extension
functor $\operatorname{Ext}_\tau$ on the category of 
noncommutative smooth
algebras, for any Fr\'echet operator ideal $\Cal K_\tau$. 
Combining
the techniques devised by Arveson and Voiculescu, we 
generalize
Voiculescu's theorem to smooth algebras and Fr\'echet 
operator ideals.
A key notion involved is $\tau$-smoothness, which is 
verified for the
algebras of smooth functions, via a noncommutative Sobolev 
lemma. The
groups $\operatorname{Ext}_\tau$ are computed for many 
examples.\endabstract
\thanks This research was supported in part by U.S. 
National Science Foundation
grant DMS9012753\endthanks
\endtopmatter

\document

\heading 1. Introduction\endheading  For a compact 
manifold $M$, the
extensions of $C(M)$ by the compact operators $\scr K(\scr 
H)$ form an abelian
group $\Ext(M)$, coinciding with the odd $K$-homology 
$K^1(M)$ \cite{BDF},
which can also be defined in terms of the elliptic 
operators \cite{A1}.  This
was a starting point of \cite K.

For a Schatten ideal $\scr L^p$, the notion of $\scr 
L^p$-smooth elements in
$\Ext(M)$ was introduced and studied in \cite{D, DV}, and 
generalized in
\cite{S1, G}.  Connes constructed the Chern character of 
extensions of smooth
algebras by $\scr L^p$ in the cyclic cohomology of the 
smooth algebras \cite C. 
When $p=1$ it recovers the trace forms in \cite{HH, CP}.

Today $\Ext$-theory for $C^*$-algebras has developed into 
a multifaced field of
fundamental importance in modern analysis, as the meeting 
ground of classical
operator theory, in particular Toeplitz operators, 
Wiener-Hopf operators, and
noncommutative differential geometry \cite C, especially 
pseudodifferential
operators, index theory, $K$-theory, and cyclic homology.

While our knowledge of topological $\Ext$-theory is rather 
complete, a natural
and fundamental problem \cite{A2, p.\ 9; D, p.\ 68; HH, 
p.\ 236} remains wide
open in noncommutative differential geometry, i.e., the 
formulation and
understanding of the extension theories of smooth 
algebras.  It would naturally
serve as the domain of the Chern character defined in 
\cite C.  Since smooth
algebras are Fr\'echet algebras, it is desirable to have 
the extension theory
constructed for Fr\'echet operator ideals.

In this note we present the analytical foundation of such 
a general theory,
which produces a functor $\Ext_\tau$ from the category of 
smooth algebras to
abelian groups, for any Fr\'echet operator ideal $\scr 
K_\tau$ \cite{P; W1,
10.11}.  Such a theory will establish a unified framework 
into which all
previous results in this direction fit.

There are two well-known technical results forming the 
cornerstone of
Brown-Douglas-Fillmore theory.  One is the celebrated 
Voiculescu's theorem
\cite{V1}, which generalizes Weyl-von Neumann-Berg theorem 
to $C^*$-algebras. 
In \cite{V2} it is extended to normed ideals and algebras 
with countable bases. 
The other is Stinespring's theorem \cite S.  We generalize 
both theorems to
Fr\'echet ideals and smooth algebras, answering a question 
in \cite D, and we
illustrate the theory in the case of smooth manifolds.

\heading 2. Smooth category and Voiculescu's 
theorem\endheading  By a {\it
smooth algebra\/} we mean a Fr\'echet $^*$-algebra 
$A^\infty$ equipped with a
norm $p$ such that $p(a^*a)=p(a)^2$ for all $a\in A^\infty$.
Often it is denoted as a pair $(A^\infty, A)$, where $A$ 
is the
$C^*$-completion of $A^\infty$ with respect to $p$.  A 
prototype is
$(C^\infty(M)$, $C(M))$ for a compact smooth manifold.

The {\it smooth category\/} is the category of all 
separable smooth algebras,
with morphisms given by $^*$-homomorphisms of Fr\'echet 
$^*$-algebras while
contractive with respect to the $C^*$-norms.

There are two key ingredients in our analysis.  Both can 
be taken for granted
for $C^*$-algebras.  One is the existence of quasi-central 
approximate
identity, identified in \cite{Ar, V2}.  The other is 
$\tau$-smoothness, which
had not appeared in literature until now.

We say a completely positive map $\phi\:A^\infty\to\scr 
L(\scr H)$ is
$\tau$-quasi-central, if there is an increasing sequence 
$\{k_m\}$ of finite
rank positive operators strongly converging to 1 such that
$$
\lim_{m\to\infty}[k_m,\phi(a)]=0\quad\text{ in }\scr 
K_\tau\text{ for every
}a\in A^\infty.
$$

Let $(A^\infty,A)$ be a smooth algebra with a dense 
$^*$-subalgebra
$A^\infty_0$, which is countably generated as a vector 
space.  A completely
positive map $\phi$ from $(A^\infty,A)$ into $\scr L(H)$ 
is $\tau$-{\it
smooth\/} (or $\scr K_\tau$-{\it smooth\/}) with respect 
to $A^\infty_0$,
if the following property holds:

If $\phi_0$ is a completely positive map from $A$ into 
$\scr L(H)$ such that
$$
(\phi-\phi_0)(A^\infty_0)\subset\scr K_\tau(H),\tag1
$$
then
$$
(\phi-\phi_0)(A^\infty)\subset\scr K_\tau(H).
$$

For $(A^\infty,A)\subset\scr L(\scr H)$, if we add
$(\phi-\phi_0)(A^\infty\cap\scr K_\tau(\scr H))\subset\scr 
K_\tau(H)$ in (1),
then we obtain the definition of a $\tau$-{\it smooth} 
$\mod\scr K_\tau$ (or
$\scr K_\tau$-{\it smooth} $\mod\scr K_\tau$) map.

A smooth operator algebra $(A^\infty, A)$ is called 
$\tau$-{\it smooth} (or
$\tau$-{\it smooth}$\mod\scr K_\tau)$ if the map 
${\roman{id}}_{A^\infty}$
is so.  The following is the generalized Voiculescu's 
theorem.

\thm{Theorem 1} Let $(A^\infty,A)$ be a separable operator 
algebra
$\tau$-smooth $\mod\scr K_\tau$ in $\scr L(\scr H)$.  Let 
$\pi$ be a
nondegenerate $\tau$-quasi-central representation of 
$(A^\infty,A)$ into $\scr
L(\scr H)$ such that $\pi|_{A\cap\scr K(\scr H)}=0$.  Then 
there are unitaries
$U_n\:\scr H\to\scr H\oplus\scr H$ such that
$$\gathered
(U^*_n(a\oplus\pi(a))U_n-a)\in\scr K_\tau(\scr 
H)\quad\text{ for every }a\in
A^\infty,\\
\lim(U^*_n(a\oplus\pi(a))U_n-a)=0\quad\text{ in }\scr 
K_\tau(\scr H)\text{ for
every }a\in A^\infty_0.\endgathered
$$\ethm

Restricting the ideal $\scr K_\tau$ and $A^\infty$ to 
various categories, one
recovers all the previous results in this direction.  The 
idea of the proof is
an extension of the techniques invented in \cite{Ar} and 
\cite{V2} for Fr\'echet
algebras.

\heading 3. Smooth extensions\endheading  A $\tau$-{\it 
smooth extension of \/}
$(A^\infty,A)$ is a pair $(\pi,P)$, where $\pi$ is a 
representation of $A$ in a
Hilbert space $\scr H$ and $P$ is a projection in $\scr H$ 
such that (i) $[\pi(
a),P]\in\scr K_\tau(\scr H)$ for $a\in A^\infty$, and (ii) 
$P\pi(A)P\cap\scr K(
\scr H)=\{0\}$.

Define a completely positive map $\phi(a)\coloneq 
P\pi(a)P,\ a\in A$.  Then
$E^\infty\coloneq\phi(A^\infty)+\scr K_\tau$ is a Fr\'echet
$^*$-algebra with the locally convex final topology 
induced by the maps
$i\:\scr K_\tau\to E^\infty$ and $\phi$, and 
$(E^\infty,E)$ is a smooth
operator algebra with a dense $^*$-subalgebra
$E^\infty_0\coloneq\phi(A^\infty_0)+\scr K_f$ countably 
generated as a vector
space.

Two $\tau$-smooth extensions are unitarily equivalent if 
the two associated
completely positive maps on $A^\infty$ are unitarily 
conjugate up to $\scr
K_\tau$-compact perturbation.  A $\tau$-smooth extension 
is degenerated if it
is unitarily equivalent to a representation.

Let $\scr Ext_\tau(A^\infty)$ be the unitary equivalence 
classes of the
$\tau$-smooth extensions of $A^\infty$.  A spatial 
isomorphism $M_2(\scr L(\scr
H))\simeq\scr L(\scr H)$ turns $\scr Ext_\tau(A^\infty)$ 
into an abelian
semigroup.  The quotient abelian monoid modulo the 
degenerate extensions will
be denoted by $\Ext_\tau(A^\infty)$.

We shall denote by $\scr Ext_{s,\tau}(A^\infty)$ those 
represented by
$\tau$-smooth completely positive maps.  Replacing $\scr 
Ext_\tau$ by $\scr
Ext_{s,\tau}$ in the above, we get the submonoid 
$\Ext_{s,\tau}(A^\infty)$ of
those consisting of $\tau$-smooth completely positive 
maps.  By a refinement of
Stinespring's theorem \cite S, we can show that 
$\Ext_\tau$ is a contravariant
functor from the category of smooth algebras to the 
category of abelian groups.

If $\tau$ is finer than $\tau'$, there is a natural 
transformation
$\alpha_{\tau,\tau'}$ from $\Ext_\tau$ to $\Ext_{\tau'}$.  
In particular, there
is a natural transformation $\alpha_\tau$ from $\Ext_\tau$ 
to the B-D-F functor
$\Ext$ for any $\tau$.  Theorem 1 implies

\thm{Corollary} If every extension in $\scr 
Ext_\tau(A^\infty)$ can be
represented by a $\tau$-smooth completely positive map and 
is
$\tau$-quasi-central if it is degenerate, then there is a 
natural isomorphism
of abelian groups
$$
\Ext_\tau(A^\infty)\simeq\scr Ext_\tau(A^\infty).\quad\bxx
$$\ethm

We check the two conditions above for the case of 
commutative smooth algebras
$(A^\infty,A)$ where $A=C(M)$, for a compact smooth 
manifold $M$ of dimension
$n$.  We may assume $M$ is embedded in $\Bbb R^N\ 
(N\leq2n$, by Whitney's
theorem).  From the deep results in \cite{V2}, we have

\thm{Theorem 2 {\rm(Voiculescu)}\rm}  Notation is as 
above.  All degenerated smooth
extensions of $(A^\infty,A)$ by $\scr L^N$ are $\scr 
L^N$-quasi-central.\ethm

For the other condition we have

\thm{Theorem 3} Let $M$ be a compact $C^{(k)}$-manifold of 
dimension $n$.
 For any Fr\'echet operator ideal $\scr K_\tau$, any
completely positive map defining extensions of $C^{(n+2+
\varepsilon)}(M)$, $C(M))$ by
$\scr K_\tau$ is always $\tau$-smooth, for any 
$\varepsilon>0$.  Here we assume
$k\geq n+2+\varepsilon$.\ethm

The key step in the proof is the following noncommutative 
Sobolev lemma.  It
implies that if a sequence $(a_{m_1,\dots,m_n})$ belongs 
to a Sobolev space
$H_s$ for sufficiently large $s$, then the quantized 
generalized function
determined by $(a_{m_{1},\dots,m_{n}})$ converges in the 
$\tau$-smooth quantized
Fr\'echet algebra.

\thm{Lemma} Let $(A^\infty,A)$ be a commutative smooth 
algebra with $n$
generators $\{x_1,\dots,x_n\}$, such that

{\rm(1)} The $C^*$-norm $\|x_j\|\leq1$, for all 
$j=1,\dots,n${\rm;}

{\rm(2)} Any element in $A^\infty$ has the form
$$
f=\sum_{(m_1,\dots,m_n)\in\Bbb 
Z^n}a_{m_{1},\dots,m_n}x_1^{m_1}\dots
x_n^{m_n},
$$
such that
$$
\sum_{(m_1,\dots,m_n)\in\Bbb
Z^n}a_{m_{1},\dots,m_{n}}(|m_1|+\dots+|m_n|)^2<\infty.
$$
For any Fr\'echet operator ideal $\scr K_\tau$, let $\phi$ 
be the completely
positive map associated to any $\tau$-smooth extension of 
$(A^\infty,A)$.  Then
$\phi$ is $\tau$-smooth.\ethm

\heading 4. Examples\endheading 

1.  Let $A=C(X)$ where $X$ is a compact, second countable, 
totally disconnected
space.  Then $A$ has a single generator.  If $A^\infty$ 
also does, then
$\Ext_{s,\tau}(A^\infty)=0$.

2.  Let $A=A^\infty=C(S^1)$.  There is a representation
$\pi\:A\longrightarrow\scr L(H)$ defining a degenerate 
$\tau$-smooth extension,
which is $\scr K$-smooth but not $\scr L^1$-smooth as a 
completely positive
map.  Thus $\Ext_{\scr L^1}$ is not a compatible functor 
for the category of
$C^*$-algebras.

3.  For any $p>1$, there is a degenerate faithful 
representation of
$(L^1(T)$,$C^*(T))$ defining a $\scr L^p$-smooth 
extension, which is not $\scr
L^p$-smooth as a completely positive map.

This shows that if $A=C^*(G)$ for a compact Lie group 
$G,L^1(G)$ is too large a
smooth subalgebra.  One needs to take, for example 
$C^\infty(G)$ instead.

4.  For any operator ideal $\scr K_\tau$, we have 
$$
\Ext_\tau C^\infty(S^1)\simeq\Ext C(S^1)\simeq\Bbb Z.
$$

Fix $\scr K_\tau\subset\scr L^p$ for some $p\geq1$.  It 
follows from \cite C
that there is a natural group homomorphism
$$
ch_m\:\Ext_\tau(A^\infty)\longrightarrow HC^{2m+
1}(A^\infty),\qquad m\geq p,
$$
such that $ch_{m+1}=S\ccirc ch_m$. Here $HC^*$ is the 
cyclic homology.

5.  Let $D$ be the unit disc.  Then $\Ext(C(D))\!=\!0$.  
However,
$\Ext_{L^1}\!C^\infty\!(D)\otimes\Bbb C$ contains as a 
direct summand the group
$HC^1(C^\infty(D))$, which is the space of all the closed 
de Rham
currents on $D$ of dimension 1.

\rem{Remarks}  Since there is no hausdorff topology on 
$\scr L/\scr K_\tau$, we
abandon the conventional formulation of $\Ext$; so a 
lifting problem does not
arise.

An attractive perspective of the smooth extension theory 
is that
$\Ext_\tau(A^\infty)$ is a differential invariant for 
appropriate $\scr
K_\tau$ and $A^\infty$.  In \cite{Km} it is shown that 
even in the B-D-F group
the class of a smooth extension may depend on the smooth 
structure.

Details will appear elsewhere (see \cite{W2, W3, W5}).  We 
plan to investigate
the topological aspect of the theory, along with the even 
degree functors in a
future work.\endrem

\heading Acknowledgment\endheading

We would like to thank J. Anderson, L. Brown, A. Connes, 
R. Douglas, J.
Horvath, J. Kaminker, and especially J. Rosenberg, and D. 
Voiculescu, for
stimulating conversations concerning various aspects of 
this work.

\enddocument